\documentclass{ws-jtca}
\usepackage{tikz}
\usetikzlibrary{patterns,patterns.meta}
\usepackage{lineno}
\usepackage[normalem]{ulem}
\begin{document}

\newcommand{\mi}{\mathrm{i}}
\newcommand{\pfrac}[2]{\frac{\partial #1 }{\partial #2 }}
\newcommand{\lr}[1]{\left( #1 \right)}
\newcommand{\ty}{  Y }
\newcommand{\ceil}[1]{\lceil #1 \rceil}
\newcommand{\intd}{\mathrm{d}}
\newcommand{\ex}{\textnormal{e}}
\newcommand{\F}{\mathbb{F}}
\newcommand{\iF}{\mathbb{F}^{-1}}

\newcommand{\BC}[1]{{\leavevmode\color{blue} Ben : #1}}
\newcommand{\MC}[1]{{\leavevmode\color{magenta} Matthew : #1}}
\newcommand{\MCo}[2]{{\leavevmode\color{red}\sout{#1} \leavevmode\color{magenta} #2}}


\markboth{M. J. King, B. T. Cox and B. E. Treeby}{Non-uniform time-stepping in k-space pseudospectral time domain models of acoustic propagation.}

%
\catchline{}{}{}{}{}
%

\title{Non-uniform time-stepping in k-space pseudospectral time domain models of acoustic propagation}

\author{Matthew J. King, B. E. Treeby and B. T. Cox}

\address{Department of Medical Physics and Biomedical Engineering, University College London,\\
London, WC1E 6BT, England.\\
\email{matthew.king@ucl.ac.uk}}

\maketitle

\begin{history}
\received{(16 June 2025)}
\end{history}

\begin{abstract}
Non-uniform time stepping in acoustic propagation models can be used to preserve accuracy or reduce computational cost for an acoustic simulation with a wave front propagating through a domain with both heterogeneous and homogenous regions, such as for a simulation of breast ultrasound tomography. The k-space correction already exist within the literature to remove numerical dispersion caused by the time stepping procedure in pseudo-spectral time domain models, but requires a uniform time step. Here we expand this correction to be able to account for a non-uniform time stepping method and illustrate the potential advantages and considerations.
A version of this Article has been submitted for review to the Journal of Theoretical and Computational Acoustics.\\

\end{abstract}

\keywords{k-space; pseudo-spectral; Time stepping.}

\section{Introduction}

When solving partial differential equations numerically a critical choice is the numerical algorithm being used. In this decision are considerations not just for the accuracy of the approach, but also the computational efficiency, memory requirements, and the aims of the output.\\

Time stepping numerical methods are often considered in acoustics for initial value problems (IVPs) when considering a heterogeneous domain due to their broadband nature and the absence of an analytic solution based on the initial data in most cases of interest. The choice of time-step size directly impacts the accuracy and computational cost. The rate of convergence of the solution in the time step size is based on the specific algorithm, while the computational cost of the problem scales linearly with the total number of time steps required for the whole simulation.\\
For the spatial domain a discrete spatial grid may be considered, with this fixed at the start of the simulation. The spatial accuracy of the solution then depends on the choice of grid size, uniform or otherwise, and the methods to compute the spatial derivatives.
In this work we consider an explicit finite difference time stepping method, coupled with a Fourier spectral method applied to the spatial derivatives though the fast Fourier Transform (FFT) allowing us to view the problem within the wave-vector $k$-space. This is applied to the first order wave equations, with time staggering between the updates for the acoustic pressure and particle velocity.\\

One advantage of using a simple finite difference method instead of a higher order scheme for the time integration is the reduction of the required storage of the problem. This is however known to introduce a phase error into the solution.
Fortunately, when considered with the spectral approach to the spatial derivatives improvements can be made. \\
To pseudo-spectral method we apply a denominator function within the finite difference scheme, these are typically considered in terms of the frequency, but transformed into a correction term in the wave-vector through the dispersion relation and then transferred onto the spatial derivative. We refer to this as the `\textit{k-space correction.}' Because the k-space correction is derived from the analytic solution within a homogenous medium the resulting solution is numerically exact both spatially and temporally to the analytic solution for the IVP. This holds true for all frequencies up to the Nyquist limit, as determined by the underlying spatial grid irrespective of the time-step size, with global stability. In heterogeneous mediums the solution is no longer exact, reintroducing some dispersion error. However, the numerical scheme including the k-space correction present lower phase error than would be observed without the inclusion of the k-space correction\cite{cox2007k}. \\

It is worth considering using non-uniform time-stepping when considering wave propagation through a heterogenous medium where any scatterers are contained within a subregion of the whole domain. This could be to reduce computational costs by introducing a larger time-step when the wave front is away from the scatterer. Alternatively, when the wave front is within the scatterer region taking a smaller time-step may be considered in order to maintain a higher accuracy. If the same simulations were run with a uniform time-step size it is understood that a trade off is required between the accuracy of the solution and computation time. An example of such a problem would be breast ultrasound tomography, though any ultrasound scans where direct contact is ill-suited would benefit from this method. Added to summarise next paragraph, moved to discussion: Similar ideas are explored through schemes with local time stepping where the spatial and temporal meshes are refined together across all time\cite{diaz2009energy}, and by considering time integration through the Runge-Kutta method where the step size is varied dynamically. In this work however we treat the temporal domain globally in space, and update our time-step size at fixed times based on the problems geometry.\\

This paper is divided into five sections. In section \ref{sect2:theory} we introduce the problem for which this paper focuses on, the first-order linearised wave equations, deriving both the k-space correction for an equal space time-stepping algorithm, but also for non-uniform time-stepping. Additionally some techniques that are required for consideration of a time-staggered grid between a system of equations. 
Section \ref{sec:exact_scheme} gives an algorithm for using the k-space corrections in this context, which is then used in the production of sections \ref{sect3:homnum} and \ref{sect4:hetnum}. \\
\ref{sect3:homnum} examines the errors from a model problem in both one-dimensional (1D) and two-dimensional (2D) homogeneous domains. 
Section \ref{sect4:hetnum} then considers the errors from heterogeneous domains with 2D examples including one drawn from photoacoustic ultrasound imaging. We discuss key points of consideration of how to change the time-stepping is specific cases and highlight the advantages of the non-uniform time-stepping over either using uniform large or small time-steps.

\section{Derivations of k-space corrections}\label{sect2:theory}

In numerical time domain models of acoustic propagation, the numerical dispersion that arises from approximating the time derivatives has been counteracted by applying corrections in the spatial Fourier (k-space) domain\cite{fornberg1978numerical,tabei2002k,mickens2020nonstandard,treeby2017nonstandard}.
Such k-space corrections have typically been defined for numerical solution schemes that use equally-spaced time-steps, ie. the time-step does not change throughout the simulation. The derivation of such a scheme will be presented in this section, to introduce the notation for Sec. \ref{sec:unequal_dt}, in which 
this derivation is extended to the case in which the time-step size can change during the simulation.\\
We consider the first-order wave equations on a staggered temporal grid. The temporal derivative is approximated with a first-order central difference scheme which is manipulated to update for the next time-step. The k-space correction is derived from the time-stepping algorithm and analytic solutions under Fourier Transform.

\subsection{Equal time-stepping}\label{sect2.1kspacederiv}

Consider the adiabatic, linearised, first-order acoustic equations in a domain with homogeneous material properties, sound speed $c_0$ and density $\rho_0$:
\begin{align}\label{FirstOrderWaves}
    \boldsymbol{U}_t(\boldsymbol{x},t)=-\frac{1}{\rho_0}\boldsymbol{\nabla}p(\boldsymbol{x},t), \qquad
    p_t(\boldsymbol{x},t)=-c_0^2 \rho_0 \boldsymbol{\nabla}\cdot\boldsymbol{U}(\boldsymbol{x},t),
\end{align}
where $\boldsymbol{U}(\boldsymbol{x},t)$ the particle velocity (the acoustic fluid velocity), and $p(\boldsymbol{x},t)$, the acoustic pressure, are functions of space $\boldsymbol{x}$ and time $t$. These can be combined to give the second order wave equation for $p$ or $\boldsymbol{U}$:
\begin{equation}\label{secondOrderU}
    \boldsymbol{U}_{tt}(\boldsymbol{x},t)=c_0^2 \nabla(\nabla\cdot \boldsymbol{U}(\boldsymbol{x},t)).
\end{equation}
When taking a pseudo-spectral time domain approach to solving equations \eqref{FirstOrderWaves} through an iterative procedure on a time staggered grid with time-step $\delta t$, we start with the central difference scheme for the time derivative,
\begin{equation} 
\frac{\boldsymbol{U}(\boldsymbol{x},t+\frac{\delta t}{2}) - \boldsymbol{U}(\boldsymbol{x},t - \frac{\delta t}{2})}{\delta t} \approx \boldsymbol{U}_t(\boldsymbol{x},t).
\end{equation}
This can then be re-arranged to update for $\boldsymbol{U}(\boldsymbol{x},t+\frac{\delta t}{2})$. This will introduce error at each time-step of order $O(\delta t^2)$. In order to remove this error a non-standard finite difference method can be used, which in the pseudo-spectral setting gives rise to the k-space correction. This is performed in the Fourier domain, by first taking the Fourier transform from $\boldsymbol{x}\to\boldsymbol{k}$. We denote the Fourier Transform as $\F$, and inverse Fourier transform as $\iF$, with $\tilde{\cdot}$ indicating the transformed variables: 
\begin{equation}
    \F\{\boldsymbol{U}(\boldsymbol{x},t)\}=\Tilde{\boldsymbol{U}}(\boldsymbol{k},t). \qquad     \F^{-1}\{\Tilde{\boldsymbol{U}}(\boldsymbol{k},t)\}=\boldsymbol{U}(\boldsymbol{x},t).
\end{equation}
The purpose of the k-space correction, $\kappa(k)$, is to satisfy
\begin{equation}\label{fixedTimeStepkappaEqn}
\frac{\Tilde{\boldsymbol{U}}(\boldsymbol{k},t+\frac{\delta t_1}{2}) - \Tilde{\boldsymbol{U}}(\boldsymbol{k},t - \frac{\delta t_1}{2})}{\delta t} = \Tilde{\boldsymbol{U}}_t(\boldsymbol{k},t) \kappa(\boldsymbol{k}),
\end{equation}
in other words, to correct for the finite different approximation. Now consider the Fourier transform of the second order linearised wave equation \eqref{secondOrderU},
\begin{equation}\label{HHEqn}
    \Tilde{\boldsymbol{U}}_{tt}(\boldsymbol{k},t)= - c_0^2 |\boldsymbol{k}|^2 \Tilde{\boldsymbol{U}} (\boldsymbol{k},t),
\end{equation}
which has the general solution
\begin{equation}\label{genWaveU}
    \Tilde{\boldsymbol{U}}(\boldsymbol{k},t)= A(\boldsymbol{k}) \ex^{ \mi c_0 \boldsymbol{k} t } + B(\boldsymbol{k})  \ex^{ -\mi c_0 \boldsymbol{k} t }.
\end{equation}
By substituting equation \eqref{genWaveU} into \eqref{fixedTimeStepkappaEqn}, the k-space correction can be seen to be
\begin{equation}\label{kappa}
    \kappa(\boldsymbol{k})= 2\frac{\sin(c_0 \boldsymbol{k} \frac{\delta t}{2})}{c_0 \boldsymbol{k} \delta t}.
\end{equation} 
Substituting \eqref{FirstOrderWaves} for the temporal derivative into \eqref{fixedTimeStepkappaEqn} gives an update equation for the particle velocity:
\begin{equation}
    \label{eq:U_equal_dt}
    \boldsymbol{U}(\boldsymbol{x},t+\tfrac{\delta t}{2}) = \boldsymbol{U}(\boldsymbol{x},t - \tfrac{\delta t}{2}) + \frac{\delta t}{\rho_0} \iF\{ \mi \kappa(\boldsymbol{k}) \boldsymbol{k} \F\{p(\boldsymbol{x},t)\} \}.
\end{equation}
The same procedure can be used to find the companion equation for the acoustic pressure:
\begin{equation}
    \label{eq:p_equal_dt}
    p(\boldsymbol{x},t+\delta t) = p(\boldsymbol{x},t) + \delta t c_0^2\rho_0 \iF\{ \mi \kappa(\boldsymbol{k}) \boldsymbol{k} \cdot \F\{\boldsymbol{U}(\boldsymbol{x},t + \tfrac{\delta t}{2})\} \}.
\end{equation}
Equations \eqref{eq:U_equal_dt} and \eqref{eq:p_equal_dt}, with the k-space correction given in \eqref{kappa}, describe a numerical scheme for solving the linearised first order acoustic wave equations \eqref{FirstOrderWaves} for a domain with homogeneous acoustic properties that is exact and unconditionally stable for any size of time-step $\delta t$, so long as that time-step is constant throughout the simulation. When the medium is not homogeneous, when the sound speed or density vary spatially, the k-space correction is no longer exact but it has been observed to reduce dispersion and so improve convergence rates\cite{treeby2012modeling,treeby2017nonstandard}.
\subsection{Non-uniform time-stepping}
\label{sec:unequal_dt}
 
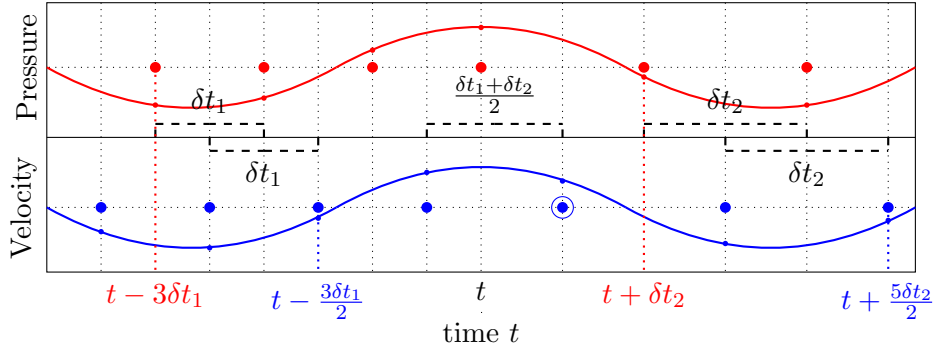
\begin{figure}[h!]
    \centering
    \begin{tikzpicture}[scale=5/7]
\draw (-8,5.5) rectangle (8,0.5);
\draw[dotted]  (8,4.3) -- (-8,4.3) node[anchor=south, rotate=90] {Pressure};
\draw (8,3) -- (-8,3);
\draw[dotted]  (8,1.7) -- (-8,1.7) node[anchor=south, rotate=90] {Velocity};
\draw[dotted]  (0,5.5) -- (0,0.5) node[anchor=north] {$t$};
\draw[dotted]  (0,-0.25) -- (0,-0.25) node[anchor=north] {time $t$};
\draw[dotted]  (-1,5.5) -- (-1,0.5) node[anchor=north] {$ $};
\draw[dotted]  (-2,5.5) -- (-2,0.5) node[anchor=north] {$ $};
\draw[dotted]  (-3,5.5) -- (-3,1.7) node[anchor=north] {$ $};
\draw[dotted]  (-4,5.5) -- (-4,0.5) node[anchor=north] {$ $};
\draw[dotted]  (-5,5.5) -- (-5,0.5) node[anchor=north] {$ $};
\draw[dotted]  (-6,5.5) -- (-6,4.3) node[anchor=north] {$ $};
\draw[dotted]  (-7,5.5) -- (-7,0.5) node[anchor=north] {$ $};
\draw[dotted]  (1.5,5.5) -- (1.5,0.5) node[anchor=north] {$ $};
\draw[dotted]  (3,5.5) -- (3,4.5) node[anchor=north] {$ $};
\draw[dotted]  (4.5,5.5) -- (4.5,0.5) node[anchor=north] {$ $};
\draw[dotted]  (6,5.5) -- (6,0.5) node[anchor=north] {$ $};
\draw[dotted]  (7.5,5.5) -- (7.5,1.7) node[anchor=north] {$ $};
\fill[red] (-6,4.3) circle (0.1cm);
\fill[red] (-4,4.3) circle (0.1cm);
\fill[red] (-2,4.3) circle (0.1cm);
\fill[red] (0,4.3) circle (0.1cm);
\fill[red] (3,4.3) circle (0.1cm);
\fill[red] (6,4.3) circle (0.1cm);
\fill[blue] (-7,1.7) circle (0.1cm);
\fill[blue] (-5,1.7) circle (0.1cm);
\fill[blue] (-3,1.7) circle (0.1cm);
\fill[blue] (-1,1.7) circle (0.1cm);
\fill[blue] (1.5,1.7) circle (0.1cm);
\draw[blue] (1.5,1.7) circle (0.2cm);
\fill[blue] (4.5,1.7) circle (0.1cm);
\fill[blue] (7.5,1.7) circle (0.1cm);
\draw[thick,blue] (-8,1.7) ..controls (-8+16/9,0.7) and (-8+32/9,0.7) .. (-8/3,1.7);
\draw[thick,blue] (-8/3,1.7) ..controls (-8/3+16/9,2.7) and (-8/3+32/9,2.7) .. (8/3,1.7);
\draw[thick,blue] (8/3,1.7) ..controls (8-32/9,0.7) and (8-16/9,0.7) .. (8,1.7);
\draw[thick,red] (-8,4.3) ..controls (-8+16/9,3.3) and (-8+32/9,3.3) .. (-8/3,4.3);
\draw[thick,red] (-8/3,4.3) ..controls (-8/3+16/9,5.3) and (-8/3+32/9,5.3) .. (8/3,4.3);
\draw[thick,red] (8/3,4.3) ..controls (8-32/9,3.3) and (8-16/9,3.3) .. (8,4.3);
\fill[red] (-6,3.6) circle (0.05cm);
\fill[red] (-4,3.73) circle (0.05cm);
\fill[red] (-2,4.62) circle (0.05cm);
\fill[red] (0,5.04) circle (0.05cm);
\fill[red] (3,4.12) circle (0.05cm);
\fill[red] (6,3.6) circle (0.05cm);
\fill[blue] (-7,1.25) circle (0.05cm);
\fill[blue] (-5,0.95) circle (0.05cm);
\fill[blue] (-3,1.5) circle (0.05cm);
\fill[blue] (-1,2.35) circle (0.05cm);
\fill[blue] (1.5,2.19) circle (0.05cm);
\fill[blue] (4.5,1.03) circle (0.05cm);
\fill[blue] (7.5,1.46) circle (0.05cm);
\draw[thick,dashed] (-6,3) -- (-6,3.25);
\draw[thick,dashed] (-4,3) -- (-4,3.25);
\draw[thick,dashed] (-6,3.25) -- (-5,3.25) node[anchor=south] {$\delta t_1$};
\draw[thick,dashed] (-4,3.25) -- (-5,3.25);
\draw[thick,dashed] (-5,3) -- (-5,2.75);
\draw[thick,dashed] (-3,3) -- (-3,2.75);
\draw[thick,dashed] (-5,2.75) -- (-4,2.75) node[anchor=north] {$\delta t_1$};
\draw[thick,dashed] (-4,2.75) -- (-3,2.75);
\draw[thick,dashed] (-1,3) -- (-1,3.25);
\draw[thick,dashed] (1.5,3) -- (1.5,3.25);
\draw[thick,dashed] (-1,3.25) -- (0.25,3.25) node[anchor=south] {$\frac{\delta t_1 + \delta t_2}{2}$};
\draw[thick,dashed] (0.25,3.25) -- (1.5,3.25);
\draw[thick,dashed] (6,3) -- (6,3.25);
\draw[thick,dashed] (3,3) -- (3,3.25);
\draw[thick,dashed] (3,3.25) -- (4.5,3.25) node[anchor=south] {$\delta t_2$};
\draw[thick,dashed] (4.5,3.25) -- (6,3.25);
\draw[thick,dashed] (4.5,3) -- (4.5,2.75);
\draw[thick,dashed] (7.5,3) -- (7.5,2.75);
\draw[thick,dashed] (4.5,2.75) -- (6,2.75) node[anchor=north] {$\delta t_2$};
\draw[thick,dashed] (6,2.75) -- (7.5,2.75);
\draw[thick,dotted,red] (3,4.3) -- (3,0.5) node[anchor=north] {$t+\delta t_2$};
\draw[thick,dotted,blue] (-3,1.7) -- (-3,0.5) node[anchor=north] {$t-\frac{3\delta t_1}{2}$};
\draw[thick,dotted,blue] (7.5,1.7) -- (7.5,0.5) node[anchor=north] {$t+\frac{5\delta t_2}{2}$};
\draw[thick,dotted,red] (-6,4.3) -- (-6,0.5) node[anchor=north] {$t-3\delta t_1$};
\end{tikzpicture}

    \caption{Illustration of a time staggered grid update for the pressure and velocity, with a change of time-step size from $\delta t_1$ to $\delta t_2$ at time $t$. Note that the effects of the change in time-stepping only impacts the update of the velocity between times $t-\frac{\delta t_1}{2}$ and $t+\frac{\delta t_2}{2}$ (circled), with all other time points having an even stencil between the point and the two previous time points, one for pressure, one for the velocity. We have chosen to illustrate an increase in the time-steps, however the method presented is equally applicable for a reduced time-step.}
    \label{fig:vairedTimeStepsIllu}
\end{figure}

Here, the derivation of the k-space correction given above will be extended to the case where the time-steps are not equal throughout the simulation but may vary, as a result of this the 'past' and 'future' in the central difference scheme are not equally spaced. An example situation, in which the time-step size changes from $\delta t_1$ to $\delta t_2$ at time $T$, is shown in Figure \ref{fig:vairedTimeStepsIllu}.
Note that the acoustic pressure and particle velocity are computed at staggered time points. This simple scenario will be used below to derive the k-space correction for non-uniform time-steps, but the correction will be applicable to the general (or perhaps extreme) case in which the time-step changes size at every step.

Considering the offset central difference algorithm approximation to the first-order derivative,
\begin{equation}
\frac{\boldsymbol{U}(\boldsymbol{x},t+\frac{\delta t_2}{2}) - \boldsymbol{U}(\boldsymbol{x},t - \frac{\delta t_1}{2})}{(\delta t_1+\delta t_2)/2} \approx  \boldsymbol{U}_t(\boldsymbol{x},t).
\end{equation}
Taking the spatial Fourier transform and once again assuming \eqref{genWaveU} as the general form for the particle velocity, it is possible to show that there does not exist a general k-space correction $\kappa(\boldsymbol{k})$ such that
\begin{equation}\label{staggeredCentralCorrection}
2\frac{\tilde{\boldsymbol{U}}(\boldsymbol{k},t+\frac{\delta t_2}{2}) - \tilde{\boldsymbol{U}}(\boldsymbol{k},t - \frac{\delta t_1}{2})}{\delta t_1+\delta t_2} =  \tilde{\boldsymbol{U}}_t(\boldsymbol{k},t)\kappa(\boldsymbol{k}).
\end{equation}
This can be seen by applying \eqref{genWaveU} into \eqref{staggeredCentralCorrection}, which gives the requirements that both
\begin{subequations}\label{kappacontradiction}
    \begin{align}
        \Bigg[A(\boldsymbol{k})\equiv0 && \text{ or } && \kappa(\boldsymbol{k})\equiv\frac{2(\ex^{\mi c_0 \boldsymbol{k} \frac{\delta t_2}{2}}-\ex^{-\mi c_0 \boldsymbol{k} \frac{\delta t_1}{2}})}{\mi c_0 k (\delta t_1 + \delta t_2)}\Bigg] & \text{ and } \\
        \Bigg[ B(\boldsymbol{k})\equiv0 && \text{ or } && \kappa(\boldsymbol{k})=\frac{2(\ex^{\mi c_0 \boldsymbol{k} \frac{\delta t_1}{2}}-\ex^{-\mi c_0 \boldsymbol{k} \frac{\delta t_2}{2}})}{\mi c_0 k (\delta t_1 + \delta t_2)}\Bigg]
\end{align}
\end{subequations}
In general, it will not be the case that $A(\boldsymbol{k})\equiv0$ or $B(\boldsymbol{k})\equiv0$, as these correspond to uni-directional propagation, so the only case for which both equations \eqref{kappacontradiction} hold true is when $\delta t_1=\delta t_2$. This reduces the problem to the scenario discussed in Section \ref{sect2.1kspacederiv}. As such it is concluded that no k-space correction exists for non-uniform time-steps when the time derivative is written in the form of \eqref{staggeredCentralCorrection}.

In order to resolve this, we note that \eqref{HHEqn} is a second-order ODE and therefore the solution space is spanned by two solutions, in this case $\tilde{\boldsymbol{U}}(\boldsymbol{k},t)$ and $\tilde{\boldsymbol{U}}_t(\boldsymbol{k},t)$, as long as both $A(k)\not\equiv0$ and $B(k)\not\equiv0$. 
With this in mind, we write the time derivative in the form
\begin{equation}\label{correctionspliteqn}
2\frac{\tilde{\boldsymbol{U}}(\boldsymbol{k},t+\frac{\delta t_2}{2}) - \tilde{\boldsymbol{U}}(\boldsymbol{k},t - \frac{\delta t_1}{2})}{\delta t_1+\delta t_2} =  \tilde{\boldsymbol{U}}_t(\boldsymbol{k},t)\kappa_1(\boldsymbol{k}) + \tilde{\boldsymbol{U}}(\boldsymbol{k},t)\kappa_2(\boldsymbol{k}),
\end{equation}
and seek expressions for the two k-space corrections $\kappa_1(\boldsymbol{k})$ and $\kappa_2(\boldsymbol{k})$. As before, consider equation \eqref{genWaveU}, now with equation \eqref{correctionspliteqn}:
    \begin{align}
        A(\boldsymbol{k}) \ex^{ \mi c_0 \boldsymbol{k} t }&\left( 2\frac{\ex^{ \mi c_0 \boldsymbol{k} \frac{\delta t_2}{2} }-\ex^{ -\mi c_0 \boldsymbol{k} \frac{\delta t_1}{2} }}{\delta t_1+\delta t_2}\right) + B(\boldsymbol{k})  \ex^{ -\mi c_0 \boldsymbol{k} t }\left( 2\frac{\ex^{ -\mi c_0 \boldsymbol{k} \frac{\delta t_2}{2} }-\ex^{ \mi c_0 \boldsymbol{k} \frac{\delta t_1}{2} }}{\delta t_1+\delta t_2}\right) \nonumber \\
        =&  A(\boldsymbol{k}) \ex^{ \mi c_0 \boldsymbol{k} t }\left(  \kappa_2(\boldsymbol{k}) + \mi c_0 \boldsymbol{k}\kappa_1(\boldsymbol{k})\right) + B(\boldsymbol{k})  \ex^{ -\mi c_0 \boldsymbol{k} t }\left( \kappa_2(\boldsymbol{k}) -\mi c_0 \boldsymbol{k}\kappa_1(\boldsymbol{k})\right).
    \end{align}
    Since this must hold true for all times $t$ we can equate the temporal coefficients to retrieve the system of equations
    \begin{subequations}
    \begin{align}
        \kappa_2(\boldsymbol{k}) + \mi c_0 \boldsymbol{k}\kappa_1(\boldsymbol{k}) =& 2\frac{\ex^{ \mi c_0 \boldsymbol{k} \frac{\delta t_2}{2} }-\ex^{ -\mi c_0 \boldsymbol{k} \frac{\delta t_1}{2} }}{\delta t_1+\delta t_2} \\
        \kappa_2(\boldsymbol{k}) - \mi c_0 \boldsymbol{k}\kappa_1(\boldsymbol{k}) =& 2\frac{\ex^{ -\mi c_0 \boldsymbol{k} \frac{\delta t_2}{2} }-\ex^{ \mi c_0 \boldsymbol{k} \frac{\delta t_1}{2} }}{\delta t_1+\delta t_2},
    \end{align}
\end{subequations}
which can be solved to give;
\begin{subequations}\label{kappa12unif}
\begin{align}
    \kappa_1(\boldsymbol{k})= 2\frac{\sin(c_0 \boldsymbol{k} \frac{\delta t_2}{2})+\sin(c_0 k \frac{\delta t_1}{2})}{c_0 \boldsymbol{k} (\delta t_1+ \delta t_2)} \\
    \kappa_2(\boldsymbol{k})= 2\frac{\cos(c_0 \boldsymbol{k} \frac{\delta t_2}{2})-\cos(c_0 \boldsymbol{k} \frac{\delta t_1}{2})}{ \delta t_1+ \delta t_2}.
\end{align}
\end{subequations}
Note that when $\delta t_1=\delta t_2$ this reduces to $\kappa_1(\boldsymbol{k})=\kappa(\boldsymbol{k})$ as in equation \eqref{kappa} and $\kappa_2(\boldsymbol{k})=0$, as expected. The time-stepping procedure becomes
\begin{align}\label{UupdateVarTimeSteps}
    \boldsymbol{U}(\boldsymbol{x},t+\tfrac{\delta t_2}{2})=& \boldsymbol{U}(\boldsymbol{x},t-\tfrac{\delta t_1}{2}) \nonumber\\
    & + \frac{\delta t_2+\delta t_1}{2} \iF \{ \tilde{\boldsymbol{U}}_t(\boldsymbol{k},t)\kappa_1(\boldsymbol{k}) +  \tilde{\boldsymbol{U}}(\boldsymbol{k},t)\kappa_2(\boldsymbol{k})\}.
\end{align}

The temptation is to think that this update equation for the particle velocity, equation \eqref{UupdateVarTimeSteps}, derived above, could be combined with the update equation for the acoustic pressure, equation \eqref{eq:p_equal_dt}, to give a numerical scheme for solving equations. However, it is not possible to use equation \eqref{UupdateVarTimeSteps} as it stands. While, from the previous two half time-steps, $\boldsymbol{U}(\boldsymbol{x},t-\tfrac{\delta t_1}{2})$ and $p(\boldsymbol{x},t)$ are known, and while equation \eqref{FirstOrderWaves} provides a means of calculating  $\tilde{\boldsymbol{U}}_t(\boldsymbol{k},t)$ from $p(\boldsymbol{x},t)$, the particle velocity at time $t$, $\tilde{\boldsymbol{U}}(\boldsymbol{k},t)$, remains unknown. One approach to overcoming this would be to make the approximation that $\tilde{\boldsymbol{U}}(\boldsymbol{k},t) \approx \tilde{\boldsymbol{U}}(\boldsymbol{k},t-\tfrac{\delta t_1}{2})$, but this would result in an inexact scheme, somewhat undermining the purpose of developing the k-space corrections. A better approach is to replace equation \eqref{correctionspliteqn} with 
\begin{equation}
\label{correctionspliteqn_2}
2\frac{\tilde{\boldsymbol{U}}(\boldsymbol{k},t+\frac{\delta t_2}{2}) - \tilde{\boldsymbol{U}}(\boldsymbol{k},t - \frac{\delta t_1}{2})}{\delta t_1+\delta t_2} =  \tilde{\boldsymbol{U}}_t(\boldsymbol{k},t)\kappa_1(\boldsymbol{k}) + \tilde{\boldsymbol{U}}(\boldsymbol{k},t - \tfrac{\delta t_1}{2})\kappa_2(\boldsymbol{k}).
\end{equation}
Note that $U_t(t)$ and $U(t-\delta t/2)$ continue to span the solution space under the condition $A(k)\not\equiv 0$ and $B(k)\not\equiv 0$. 
Proceeding exactly as in equations \eqref{correctionspliteqn}-\eqref{kappa12unif} gives the k-space corrections as
\begin{subequations}\label{kappa12stag}
\begin{align}
        \kappa_1(\boldsymbol{k})=&  \frac{2}{ \delta t_1+ \delta t_2}\left( \frac{\sin(c_0 \boldsymbol{k} \frac{\delta t_2}{2})}{c_0 \boldsymbol{k}}+\frac{\sin(c_0 \boldsymbol{k} \frac{\delta t_1}{2})}{c_0 \boldsymbol{k}}\frac{\cos(c_0\boldsymbol{k}\frac{\delta t_2}{2} )}{\cos(c_0\boldsymbol{k}\frac{\delta t_1}{2} )} \right) \\
        \kappa_2(\boldsymbol{k})=& \frac{2}{ \delta t_1+ \delta t_2}\left(\frac{\cos(c_0 \boldsymbol{k} \frac{\delta t_2}{2})}{\cos(c_0 \boldsymbol{k} \frac{\delta t_1}{2})} -1 \right),
\end{align}
\end{subequations}
with the particle velocity update equation now given by
\begin{align}\label{eq:U_update2}
    \boldsymbol{U}(\boldsymbol{x},t+\tfrac{\delta t_2}{2})= \boldsymbol{U}(\boldsymbol{x},t-\tfrac{\delta t_1}{2})
    + \frac{\delta t_2+\delta t_1}{2} \iF \{ \tilde{\boldsymbol{U}}_t(\boldsymbol{k},t)\kappa_1(\boldsymbol{k}) +  \tilde{\boldsymbol{U}}(\boldsymbol{k},t-\tfrac{\delta t_1}{2})\kappa_2(\boldsymbol{k})\}.
\end{align}

\section{Exact scheme for non-uniform time-steps}
\label{sec:exact_scheme}
In the section above, the k-space corrections were derived for the case illustrated in Fig.\ \ref{fig:vairedTimeStepsIllu} in which the time-step changed at time $t$ from $\delta t_1$ to $\delta t_2$. As the update equations require only the pressure and particle velocity at the previous time-step (and half time-step), the size of the time-step could change on every iteration. In other words, the time-steps could all be different: $\delta t_1, \delta t_2,...,\delta t_{n-1}, \delta t_n,...$. Using equations \eqref{FirstOrderWaves}, \eqref{eq:U_update2}, and \eqref{eq:p_equal_dt}, the update equations for an exact scheme (when the medium is homogeneous) can be written as:
\begin{equation}\label{eq:U_update_dtn}
    \boldsymbol{U}(\boldsymbol{x},t+\tfrac{\delta t_n}{2})= \boldsymbol{U}(\boldsymbol{x},t-\tfrac{\delta t_{n-1}}{2})
    + \frac{\delta t_n+\delta t_{n-1}}{2} \iF \{     
    - \kappa_1(\boldsymbol{k})\frac{\mi \boldsymbol{k}}{\rho_0}\tilde{p}(\boldsymbol{k},t)
    + \kappa_2(\boldsymbol{k})\tilde{\boldsymbol{U}}(\boldsymbol{k},t-\tfrac{\delta t_{n-1}}{2})\}.
\end{equation}
\begin{equation}
    \label{eq:p_equal_dtn}
    p(\boldsymbol{x},t+\delta t_n) = p(\boldsymbol{x},t) + \delta t_n c_0^2\rho_0 \iF\{ \mi \kappa(\delta t_n,\boldsymbol{k}) \boldsymbol{k} \cdot \tilde{\boldsymbol{U}}(\boldsymbol{k},t + \tfrac{\delta t_n}{2})\},
\end{equation}
with $\kappa_1$ and $\kappa_2$ given by equations \eqref{kappa12stag}, and $\kappa$ by equation \eqref{kappa}. It is this algorithm, along with initial conditions, that has been implemented to provide the numerical results presented in sections \ref{sect3:homnum} and \ref{sect4:hetnum}.
To implement initial conditions, both the pressure and the particle velocity are given at the same initial time $t=0$, we continue to be able to make use of the equation \eqref{eq:U_update_dtn} setting $\delta t_0=0$ in the computation of $U(\boldsymbol{x},\frac{\delta t_1}{2})$. Setting $\delta t_N=0$ at the end of the simulation also allows for the retrieval of the final time particle velocity through \eqref{eq:U_update_dtn}.
With fixed time-stepping, in order to implement the initial condition on the staggered grid we are instead required to implement a negative half time-step for the velocity as an adjusted initial condition. Typically this has been performed by restricting the initial velocity to be $U(\boldsymbol{x},0)=0$ and setting the negative half time-step to be the negative of the positive half time-step $U(\boldsymbol{x},-\delta t/2)=-U(\boldsymbol{x},\delta t/2)$, and adjusting the equivalent to equation \eqref{eq:U_update_dtn} accordingly\cite{} \textbf{Which k-wave paper}.
More generally however this could also be done by setting 
\begin{equation}   
\label{eq:exact_initial_conditions}
\boldsymbol{U}(\boldsymbol{x},-\delta t_0/2) = \iF\{ \cos(ck\delta t_0/2) \F\{ \boldsymbol{U}_0(\boldsymbol{x}) \} + \frac{\sin(ck\delta t_0/2)}{ck} \left(\frac{-\mi \boldsymbol{k}}{\rho_0} \right)\F\{ p_0(\boldsymbol{x}) \}\}.
\end{equation}
which can be found by considering the backwards half time-step using equation \eqref{eq:U_update_dtn} and the appropriate substitutions for $\kappa_1$ and $\kappa_2$.

\section{Homogeneous medium example: exact scheme}\label{sect3:homnum}

\begin{figure}[h!t]
    \centering
    \includegraphics[width=\textwidth]{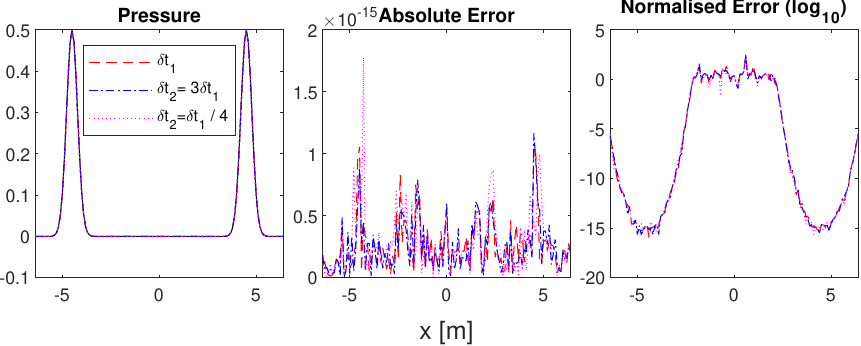}
    \caption{(Left) Pressure distribution in a Homogeneous domain after time $T$. Plotted for the analytic solution (black solid line), for a single time-step size (red dashed line), with an increased time-step at time $T/3$ (blue dot-dashed line) and with a decreased time-step at time $T/3$ (magenta dotted line). Also plotted; (Middle) The absolute difference between the analytic solution and each time-stepping approach. (Right) The normalised error between the analytic solution and each time-stepping approach under $\log_{10}$ scaling.}
    \label{Fig:1Dhomog}
\end{figure}

In this section we show that the algorithms described in Section \ref{sec:exact_scheme} produce the correct solution to equations \eqref{FirstOrderWaves} to machine precision. Additionally we illustrate that not using equation \eqref{eq:U_update_dtn} results in additional errors. We consider a 1D and a 2D homogeneous domain with a an initial pressure following a bell curve $\ex^{-(\boldsymbol{x}/4\delta x)^2}$ and $\ex^{-(\boldsymbol{x}^2 + \boldsymbol{y}^2)/{(4\delta x)^2}}$ respectively. This is performed in a finite 1D [12.9m, $\delta x=0.1$m] and 2D [$12.9\times 12.9 $m$^2$, $\delta x=\delta y = 0.1$m] domain with a perfectly matched layer (PML)\cite{treeby2011k} around the outside of the domain. All computations have been performed under double precision arithmetic, with the initial pressure being effectively compactly supported as a result.\\
We have compared the numerical results for the propagating wave at a final time $T$ when solved with uniform time-stepping, with $\kappa$ given by \eqref{kappa}, as well as with a time-step that is changed at time $T/3$ with both an increased $\delta t_2= 3 \delta t_1$ and decreased time-step $\delta t_2=\delta t_1 /4$. These are additionally compared in 1D against the analytic solution given by;
\begin{equation}
    p(\boldsymbol{x},t) = \frac{1}{2}\left(p^{(0)}(x+c_0 t) + p^{(0)}(x-c_0 t)\right) 
\end{equation}
We consider a time-step of $5ms$ over 900 time-steps to $T=4.5s$.
Figure \ref{Fig:1Dhomog} shows the pressure distribution at the final time for the 1D problem, plotting the exact solution and the solution generated with 900 equal time-steps, as well as with the adjusted time-stepping at time $2T/3$. In addition to this we additionally plot the error from the three cases and the analytic solution at the final time.
Comparing the two types of errors, the absolute difference and relative errors, 
it is important to note that when the relative error is at its largest (magnitude of 1-100)  this is when the absolute difference is at its smallest, due to analytic solution being near or equal to 0 in this region, in contrast, in the regions where the pressure is largest, which retrieves an increased error, is exactly where the relative error is at its smallest.
While each time-stepping algorithm performs similarly it is actually the method with the smallest time-steps that performs worst, with consistently the largest error. It should be noted that this error is still within machine precision of double precision arithmetic and occurs due to the accumulation of small errors from the FFT and its inverse, which is increased as a smaller time-step requires more iterations to reach the same final time.

\begin{figure}[h!t]
    \centering
    \includegraphics[width=\textwidth]{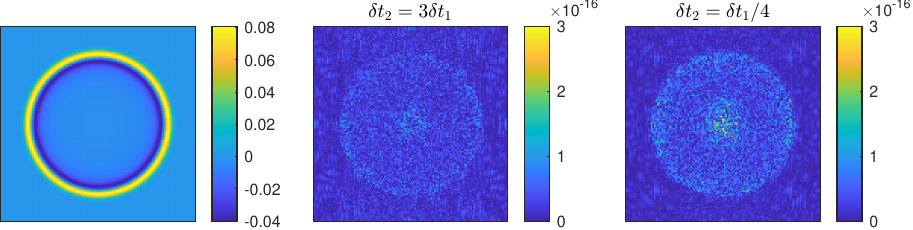}
    \caption{(Left) Pressure distribution in a homogeneous domain after time $T$ with a single size of time-step $\delta t=\delta t_1=5ms$. When the time-step is changed  this change occurs at time $2T/3$. This is for a $\delta t_2= 3\delta t_1$. (Middle) and with $\delta t_2= 3\delta t_1$ (Right.) Both the (Middle) and (Right) plots show the absolute error compared to using a single time-step at time $T$.}
    \label{Fig:2Dhomog}
\end{figure}

Figure \ref{Fig:2Dhomog} compares the single time-step to both an increased and decreased $\delta t_2$ now in a 2D domain with a time-step of $5ms$ over 900 time-steps to $T=4.5s$.

Figure \ref{Fig:2Dhomog} observes similar behaviour as in Figure \ref{Fig:1Dhomog}, with a larger error when considering smaller time-steps, due to the summation of small errors. In both cases the errors are still within double precision, and are largely contained within the region contained by the wave front. 

\begin{figure}[h!t]
    \centering
    \includegraphics[width=\textwidth]{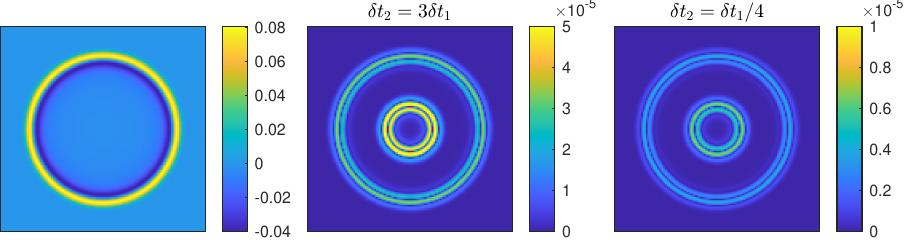}
    \caption{As Figure \ref{Fig:2Dhomog}, with $\kappa_1(\boldsymbol{k})=\kappa(\boldsymbol{k};\delta t_2)$ with the change in time-step not corrected for.}
    \label{Fig:2DhomogNoCorr}
\end{figure}

In comparison, Figure \ref{Fig:2DhomogNoCorr} examines the effects of not including the term $\kappa_2(\boldsymbol{k})$ and instead just swapping from $\kappa(\boldsymbol{k};\delta t_1)$ to $\kappa(\boldsymbol{k};\delta t_2)$. 

While it may be observed in the left plot that the pressure distributions are visually similar to that observed for the correct k-space correction for non-uniform time-steps, observing the errors between these solutions and the single time-step solution reveals the presence of a new reflected wave. At the time of the change of time-step $(2T/3)$ on top of the original outgoing wave as would be expected, an inward travelling wave is produced which alters the amplitude of the wave front such that the maximum error increases as the wave reaches the middle of the domain. Increasing either of the time-step sizes increases the resulting errors, with particularly large time-steps resulting in an error comparable to the initial pressure distribution. It should be noted that unlike in the case of the correct k-space correction taking small time-steps now increases the accuracy, though the reflected wave is still produced. 

A final approach to accounting for the correction with minimal additional computational cost would be to only change the form of $\kappa_1$ in line with \eqref{kappa12stag}, setting $\kappa_2=0$. Doing so it is found that similar errors to Figure \ref{Fig:2DhomogNoCorr}, however is only observed in half of the domain, either only accounting for the left or right propagating waves respectively, since this would be equivalent to setting $\kappa_1$ according to \eqref{kappacontradiction}, where the other condition requires $A(\boldsymbol{k})=0$ or $B(\boldsymbol{k})$=0.

\section{Heterogeneous medium example\label{sect4:hetnum}}

While the k-space correction $\kappa$ is exact for an acoustically homogeneous medium, this is no longer the case when the medium is heterogeneous (which is often required in practical applications).
When considering a stepped domain, i.e the sound speed $c_0$ and mean density $\rho_0$ are constants within disjoint subdomains, away from the heterogeneities the equations reduce to that of the homogeneous problem. As a result of this the resulting solution with the k-space correction will no longer be exact. We continue to use the k-space correction unchanged however as it continues to reduce the dispersion error and converges as the time-step size is taken to 0.

It is of note however that the k-space correction depends on the sound speed under the Fourier transform. To simplify this further and avoid the incorporation of the convolution operator a fixed reference sound speed is chosen and used across the whole domain for computation of $\kappa$. Analysis and discussion of the dispersion errors and the choice of the reference sound speed exist within the literature\cite{treeby2017nonstandard}. When considering changing the time-step, it is important to choose to do this when the wave peak is located within a region with sound speed equal to the reference used for the k-space correction, this prevents the introduction of additional dispersion errors that would be propagated into the domain.\\

\begin{figure}
    \centering
    \includegraphics[width=\textwidth]{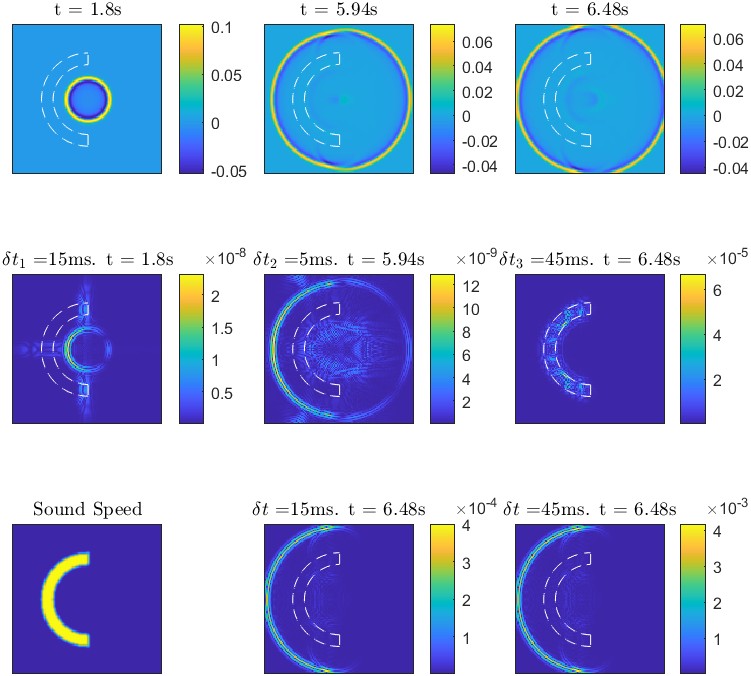}
    \caption{(Middle left) Schematic of the sound speed (equivalent to the mean density) for the simulation. (Top) Pressure at indicated time, simulated with time-steps of size $\delta t=5ms$. These are given as a reference for the other plots. (Middle middle and right) The absolute difference between the pressure perturbations at the given times from simulations run with the single time-step $15ms$ and $45ms$. (Bottom) The absolute error at the indicated times for a simulation run with $\delta t_1=15ms$, $\delta t_2=5ms$ and $\delta t_3=45ms$ with time-step changes at 1.8s and 5.94s. Also indicated (dashed lines) is the boundary between the regions of the domain with different sound speeds and density.}
    \label{Fig:2DHet}
\end{figure}

Considering a 2D domain that contains a lower sound speed within a half ring centred at the origin, and a mean density chosen such that $c_0^2 \rho_0$ is constant everywhere, for an initial pressure as used for figures \ref{Fig:2Dhomog} and \ref{Fig:2DhomogNoCorr}. Figure \ref{Fig:2DHet} is generated from simulations run with time-step sizes $\delta t=5ms$, $\delta t=15ms$ and $\delta t=45ms$. The pressure field (Figure \ref{Fig:2DHet} top row), and the comparison errors for using larger time-steps (Figure \ref{Fig:2DHet} bottom row, middle and right) have been plotted. Finally (Figure \ref{Fig:2DHet} middle row) we have considered the simulation run with non-uniform time-steps, running the simulation with $\delta t_1=15ms$ until time $1.8s$, then changing the time-step to $\delta t_2=5ms$, until time $5.94s$, and finally $\delta t_3$ until the final time at $6.48s$, plotting the error at each time. This example was constructed such that the time-step changes occurred when the wave front is outside of the half ring, resulting in the smallest time-steps being used as the wavefront crosses the heterogeneities. In addition to this the change of time-stepping occurs when the wave front is contained within the region with sound speed matching the reference sound speed.
Immediately it can be observed that the non-uniform time-stepping solution performs better than either of the larger time-step solutions, despite making use of the same time-steps at different times, in particular at the point of changing from $\delta t_1$ to $\delta t_2$ and then $\delta t_2$ to $\delta t_3$ at $1.8s$ and $5.94s$ respectively the difference between this solution and the solution taking the single small time-step $\delta t=\delta t_2=5ms$ would be negligible under single precision. This continues to be the case for the outward propagating wave, which is not visible on the plotted scale at time $6.48s$ where the simulation ends. At this final time a significant increase in the error is observed, though it is contained within the half ring. This is due to the change in time-step being incorrectly accounted for in the sound speed region with lower (than the reference) sound speed. This results in the reflected waves gaining a small amount of dispersion error. It should be noted however that this new error is 1/100th the size of the errors observed for using larger time-steps alone.\\

\begin{figure}
    \centering
    \includegraphics[width=\textwidth]{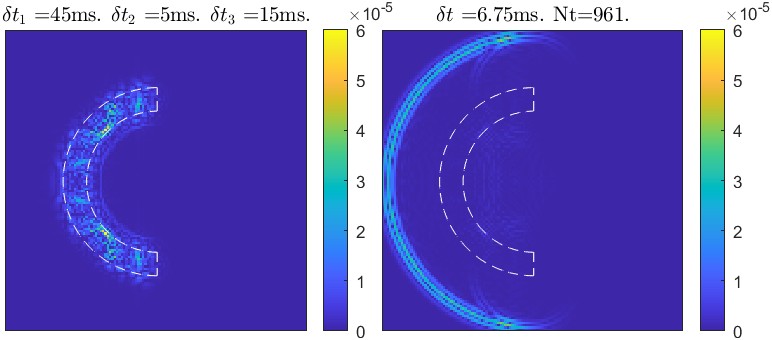}
    \caption{ Absolute errors for simulations run as for Figure \ref{Fig:2DHet} comparing the simulation with with time-steps $5ms$ to the simulation run with a non-uniform time-step $\delta t_1=15ms$, $\delta t_2=5md$ and $\delta t_3=45ms$, and for uniform time-steps $\delta t=6.75ms$ such that both simulations have the same total number of time-steps.}
    \label{Fig:2DHetAve}
\end{figure}

To give further comparison between the balance of numerical accuracy and computational cost Figure \ref{Fig:2DHetAve} plots the same simulation as given for Figure \ref{Fig:2DHet}, but with both the non-uniform time-stepping and with uniform time-stepping, set such that both simulations run for the same number of total time-steps. This ensures the improved error is not just due to the increased number of time-steps. \\
Firstly, it should be noted that both simulations do produce the same magnitude of error when compared to the simulation run with small time-steps, and overall the simulation run with uniform time-steps actually produces a lower maximum error. However, this error occurs across the entire wave front, unlike when using the non-uniform time-steps which contains the error within the heterogeneous region of the domain. As a result of this, while the accuracy overall may be considered as marginally better, if a sensor were to be placed outside of the half ring the resulting detected signal would observe larger errors than if produced for the non-uniform time-stepping, further highlighting how this method can increase accuracy even when preserving computational cost.\\

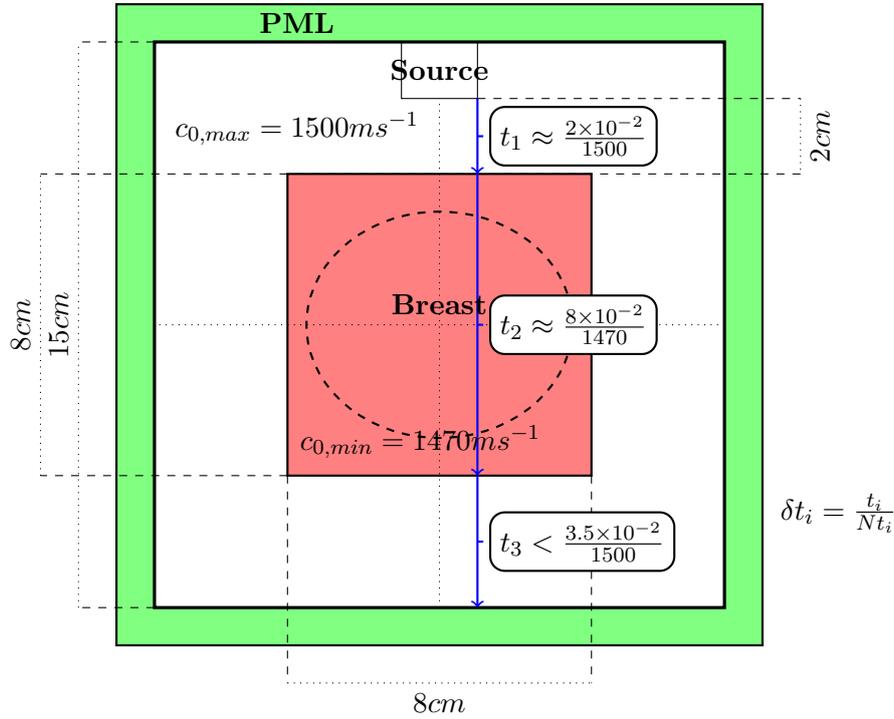
\begin{figure}
    \centering
    \begin{tikzpicture}[scale=1/2]
		\filldraw[fill=green!50, thick] (-8.5,-8.5) rectangle (8.5,8.5);
		\filldraw[fill=white, very thick] (-7.5,-7.5) rectangle (7.5,7.5);
		\filldraw[fill=red!50, thick] (-4,-4) rectangle (4,4);
		\draw[] (-1,7.5) rectangle (1,6);
		\draw[thick,dashed] (0,0) ellipse (3.5 and 3) ;
		\draw[dotted] (-9.5,-7.5) -- node[above, rotate=90, midway] {$15cm$} (-9.5,7.5);
		\draw[dotted] (-10.5,-4) -- node[above, rotate=90, midway] {$8cm$} (-10.5,4);
		\draw[dotted] (-4,-9.5) -- node[below, midway] {$8cm$} (4,-9.5);
		\draw[dotted] (9.5,6) -- node[below, rotate=90, midway] {$2cm$} (9.5,4);
		\draw[dashed] (-7.5,7.5) -- (-9.5,7.5);
		\draw[dashed] (-7.5,-7.5) -- (-9.5,-7.5);
		\draw[dashed] (9.5,4) -- (-10.5,4);
		\draw[dashed] (-4,-4) -- (-10.5,-4);
		\draw[dashed] (1,6) -- (9.5,6);
		\draw[dashed] (4,-4) -- (4,-9.5);
		\draw[dashed] (-4,-4) -- (-4,-9.5);
		\draw[dotted]  (7.5,0) -- (-7.5,0);
		\draw[dotted]  (0,-7.5) -- (0,6);
		\draw[thick, blue] (1,5) -- +(0.4em,0);
		\draw[thick, blue] (1,0) -- +(0.4em,0);
		\draw[thick, blue] (1,-5.75) -- +(0.4em,0);
		\draw[thick, blue, ->] (1,6) -- node[draw, rounded corners = 0.5em, right = 0.4em, midway, color = black, fill = white] {${t_1\approx\frac{2\times10^{-2}}{1500}}$} (1,4);
		\draw[thick, blue, ->] (1,4) -- node[draw, rounded corners = 0.5em, right = 0.4em, midway, color = black, fill = white] {${t_2\approx\frac{8\times10^{-2}}{1470}}$} (1,-4);
		\draw[thick, blue, ->] (1,-4) -- node[draw, rounded corners = 0.5em, right = 0.4em, midway, color = black, fill = white] {${t_3<\frac{3.5\times10^{-2}}{1500}}$} (1,-7.5);
		\node at (0,6.75) {\textbf{Source}};
		\node at (-3.75,8) {\textbf{PML}};
		\node at (0,0) [above] {\textbf{Breast}};
		\node at (-3.75,4.5) [above] {$c_{0,max}=1500ms^{-1}$};
		\node at (-0.5,-3.8) [above] {$c_{0,min}=1470ms^{-1}$};
		\node at (10.5,-5) {${\delta t_i = \frac{t_i}{Nt_i}}$};
\end{tikzpicture}
    \caption{Basic Schematic for the computation of the values of $\delta t$ as used for the generation of Figure \ref{Fig:2DHetBioMed} based on distance between the the source and phantom, the width of the phantom, and the sound speeds. The region containing the phantom (red) is considered with the minimum sound speed, and $Nt_2$ will be chosen to be large to maintain accuracy through this region, while the rest of the domain (white) considers the largest possible sound speed in the region and will have $Nt_1$ and $Nt_3$ chosen to reduce computational cost.}
    \label{fig:breastSchematic}
\end{figure}

\begin{figure}
    \centering
    \includegraphics[width=\textwidth]{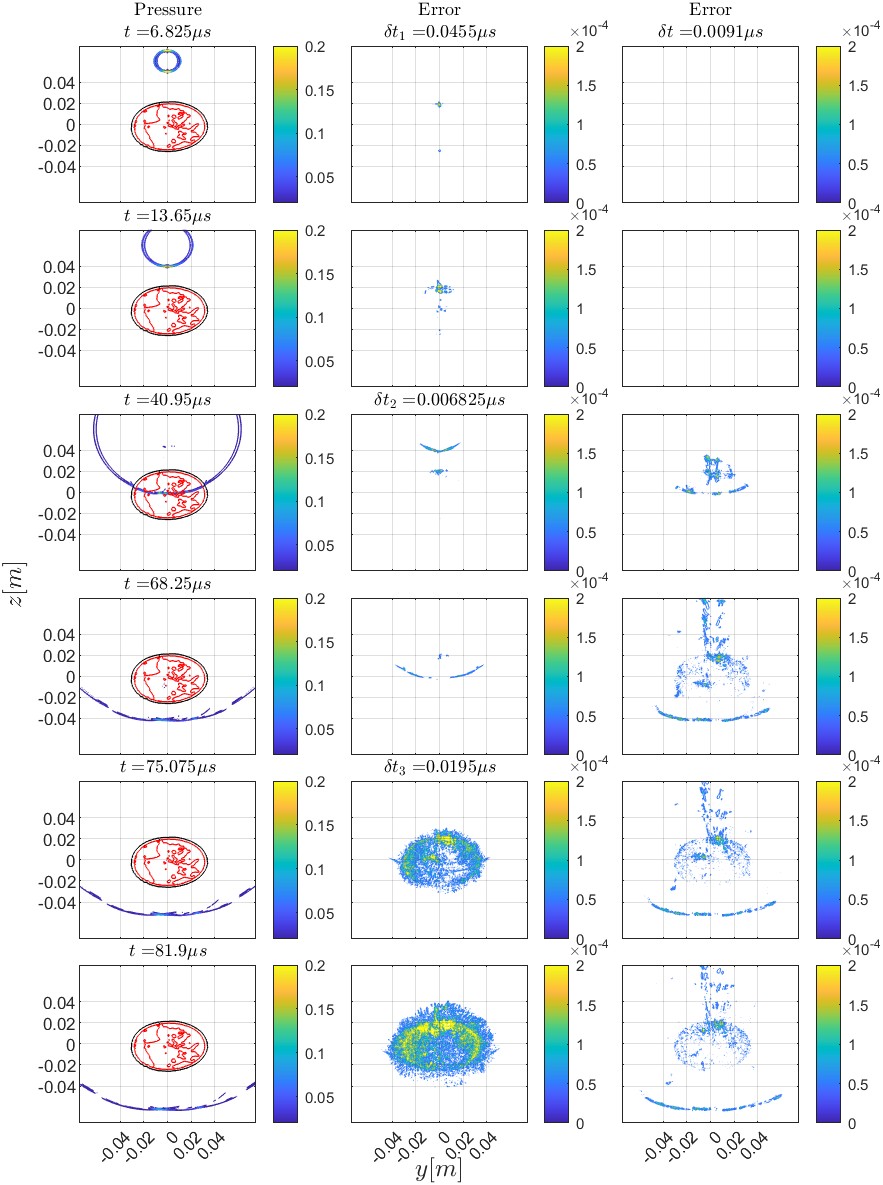}
    \caption{(Left) Reference pressure perturbation at each time computed with uniform time-steps $\delta t=7.5$ns. (Middle, Right) Absolute error compared to the reference pressure perturbation for the time-stepped solution (Middle, $\delta t_1=30$ns, $\delta t_2=7.5$ns and $\delta t_3=15$ns) and the solution with uniform time larger time-steps (Right, $\delta t=10$ns). Profiles given across central slice. Also indicated are the locations of the different tissue types in the breast phantom (Left).}
    \label{Fig:2DHetBioMed}
\end{figure}

In the simulations used to produce both Figures \ref{Fig:2DHet} and \ref{Fig:2DHetAve} the choice to change the time-step size within the region of higher sound speed was not arbitrary. If instead the simulation was arranged such that the change in time-stepping occurred in the region of lower sound speed (not equal to the reference sound speed), this results in errors equivalent to those of the original time-stepping produce for the equivalent time-step size, suggesting caution must be taken in choosing when to change the time-step, and possibly the choice of reference sound speed, in conjunction with pre-existing literature\cite{treeby2017nonstandard}. In particular smaller time-steps are more desirable on regions of the domains that exhibit the heterogeneities, and so choosing the reference sound speed based on the layout of the domain may be applicable as long as stability is accounted for. When constructing the numerical simulation described by Figure \ref{fig:breastSchematic} we have chosen $c_{0,ref}$ based on the location of the wavefront when we change time-steps ($1500ms^{-1}$). This simulation was used for the generation of Figure \ref{Fig:2DHetBioMed}.\\

Figure \ref{Fig:2DHetBioMed}, as depicted in Figure
\ref{fig:breastSchematic} considers a simulation run over a $2D$ computational domain with a spatial grid of size $150mm\times 150mm$ with a grid spacing of $0.2mm$ in each direction. A 2D slice of a 3D breast phantom\cite{lou2017generation} surrounded by water is simulated within the domain and the sound speed and mean densities are given accordingly. A 2D slice of the 3D phantom has been used in order to reduce computational complexity and simulation run time though the process for scaling the problem to 3D requires no additional analytic work. Note that acoustic absorption has not been included in this simulation.\\

An initial pressure is produced at the edge of the domain and propagated through the phantom with the pressure stored at times $6.825\mu$s $13.65\mu$s, $40.95\mu$s, $68.25\mu$s, $75.075\mu$s and $81.9\mu$s. The numerical setup emulates a photoacoustic ultrasound breast scan where ultrasound transducers are placed in a ring around the breast, with each producing an acoustic pulse in turn while the rest act as a sensors\cite{qu2019study}. In these scans the ultrasound wave propagates through water before scattering off the various types of breast tissues, before continuing to propagate through water to the sensors. Setting the reference sound speed to that of the water and using larger time-steps between times $(0,13.65 \mu s)$ and $(68.25\mu s,81.9 \mu s)$, with these times calculated according to the sound speed, and the location of the phantom, and smaller time-steps when the wave is travelling through the phantom results in a simulation that reduces the number of time-steps the solution is propagated over while maintaining a higher accuracy across the scattering region for the wave front.\\ 
In the production of Figure \ref{Fig:2DHetBioMed}, we have used smaller time-steps of uniform $6.825$ns for the reference solution (Figure \ref{Fig:2DHetBioMed}, left column). Varied between the time-step sizes of $45.5$ns, $6.825$ns and $19.5$ns (Figure \ref{Fig:2DHetBioMed}, middle column). And run the simulation with uniform intermediate time-steps of $9.1ns$ (Figure \ref{Fig:2DHetBioMed}, right column). This results in both the simulation with non-uniform time-steps and uniform intermediate time-steps  to have the same number of total iterations between the start and final time. The process for choosing initial values of the number of time-steps and time-step size is described by Figure \ref{fig:breastSchematic}.

When using smaller time-steps a total of 12000 iterations were required, while using the non-uniform time-steps as described above requires a reduced 9000 iterations, with the intermediate uniform time-step size also using 9000 time-steps. In cases where a larger scanner could be used, for example the PAM3 system\cite{dantuma2023fully} uses a hemisphere with a diameter of 26cm, the reduction in the number of time-steps will grow, additionally we have not continued the simulation to a point where the wave front has propagated through the entire domain which would further increase the computational saving by using larger time-steps. The choice of times where we have changed the time-step size also plays a role in the computational saving, with the first change occurring at $13.65\mu s$ (Figure \ref{Fig:2DHetBioMed} row 2) with the wave front significantly distanced from the phantom, allowing potential further computational savings. \\

Comparing between uniform intermediate time-steps and non-uniform time-steps, the errors are comparable in magnitude. It can however be observed once again (see Figure \ref{Fig:2DHetAve}) that in the case of the non-uniform time-steps the errors are contained away from the wave front, primarily inside the scatterer, but are distributed across the domain for the uniform time-step case, with the largest errors occurring on the wave front. \\

It should be noted that it while it may have been anticipated that that the error would have have been to machine precision, until the wave front reaches the breast phantom, small errors are observed at times $6.825\mu$s and $13.65\mu$s. This is due to the inhomogeneous sound speed in the domain along with the global computation of the spatial derivatives. This impacts the accuracy through the required choice of a reference sound speed for the computation of $\kappa_1$ and $\kappa_2$. The reference sound speed was chosen as the sound speed in water (1500ms$^{-1}$) in order to minimise the error due to the change of time-stepping.

\section{Concluding Remarks}
In this paper we have derived a k-space correction for the evaluation of an offset central difference finite difference operator for the first order wave equations. This correction is used when the time-stepping is varied between the past and the future, and consist of two terms. The first of these acts as a generalisation of the k-space correction already in the literature, valid for constant time-stepping and multiplies the derivative at the current time. The second term is zero when the time-step is constant, and multiplies the variable at the current time, or can be further modified to act on the variable at the past evaluation time. Varying the time-step allows for the computational cost to be reduced over the whole domain by taking larger time-steps without sacrificing accuracy in regions where heterogeneity's may occur in both the sound speed and the mean density where small time-steps may be taken instead. 

The new correction term is only applied at the time the time-stepping is varied, and the approach with analytic derivation is given for both when staggered time-steps \eqref{kappa12stag}, or uniform time-steps \eqref{kappa12unif} are used. The resulting time-stepping algorithms remains exact in a homogeneous domain as illustrated in Figures \ref{Fig:1Dhomog} and \ref{Fig:2Dhomog} for a time staggered grid for both a 1D and 2D domain.


An application of non-uniform time-stepping includes domains with heterogenous regions, particularly in the mean density and sound speed. In a 2D domain containing a half ring of varied mean density and sound speed the method has been observed to maintain accuracy compared to just using larger time-steps, particularly at the wave front, with errors being trapped within the heterogeneous regions and not propagating into the far field. When the sound speed is varied it is know that a reference sound speed must be chosen for the computation of the k-space correction, and this remains true for the non-uniform time-stepping. In particular it is important that the time-step change occurs when the wave front is contained within the region with sound speed matching the reference sound speed for minimal error.

Varying the time-step in order to save computational cost, or increase accuracy is not new, and is considered in the implementation of local time stepping, and Runge-Kutta methods. When using local time stepping both the underlying spatial and temporal mesh are taken to be significantly courser within a region of the domain in order to increase accuracy of the simulation, as a result the computational cost is reserved for the regions of the domain with increased sensitivity such as through heterogeneities. The distinction to the approach we present here is that we are constructing a numerical scheme where the time-step size can be changed across the whole domain in order to continue the simulation, as opposed to within a specific region of the domain for all time.
More similar to our approach would be approaching the time integration through the Runge-Kutta method, where the step step size varied dynamically based on the on variation of the function between earlier time points, however here (Section \ref{sect4:hetnum}) we are considering scenarios where the time step size is instead adapted only at fixed points in time based on the problems geometry and the known sound speeds. A method of dynamically adapting the time stepping could be considered in future work, and the numerical scheme given here would remain accurate.\\

Highlighting the applicability of this method to applications a 2D simulation has been run in the generation of Figure \ref{Fig:2DHetBioMed} where a pulse is sent through a breast phantom within water, reminiscent of photoacoustic imaging techniques, in which a ring of transducers pulse waves through the breast in sequence to reconstruct an image of the breast tissues. With only a small amount of information regarding the breast diameter and location as well as expected sound speeds times could be chosen for the simulation to reduce the number of time-steps required under non-uniform time-stepping. Comparing the resulting errors to a constant time-stepping with the same number of total time-steps it is once again observed that varying the time-steps can result in a maintained accuracy at the wave front. Additionally, the error between the two simulations remain comparable with the non-uniform time-stepping displaying a larger maximum error, and the uniform time-stepping observing a higher distribution of the error across the domain. \\

In order to extend this work further considerations for the derivation of the k-space correction terms under the temporal Fourier transform, with more general dispersion relations, would generalise the results here to other problems. Additionally, extending the problem for biomedical application, power law absorption would be added and the effects of this on the resulting error compared to the error from a single time-step size would wish to be considered. It is noted that the dispersion relation used for the derivation of the k-space correction terms would no longer be exact. Alongside this a more in-depth study of the choice of the reference sound speed considering the change of time-step and the location of the wave front in detail.\\

In a similar approach to that used here the phase error in the finite difference time stepping algorithm, which depends on the underlying frequencies, can be corrected for using a non-standard finite-difference scheme\cite{mickens2020nonstandard}. This is done by introducing a denominator function which corrects for a single frequency for a specific time-stepping procedure in a homogeneous domain, the result is highly accurate for single-frequency or narrowband problems. This method alone is not be valid here due to the presence of broadband frequencies, however connects directly to the k-space correction, through the relationship between the frequency and the wave-vector through the dispersion relation\cite{treeby2017nonstandard}. This relationship is further explained in Appendix \ref{AppA}. As such further investigation of non-standard finite difference schemes through the study of 'exact' or 'best' finite difference models may reveal additional improvements as the current model is extended to include non-linear terms and power-law absorption.

\section*{Acknowledgements}
This research was supported by the Engineering and Physical Sciences Research Council (EPSRC), UK, under grants EP/W029324/1 and EP/T014369/1.

\appendix
 \section{Non-standard Finite Difference Methods}\label{AppA}
 
 A non-standard first order finite difference method, as described by Mickens\cite{mickens2020nonstandard}, are approximations of the derivative of, as an example, the form
 \begin{subequations}
 \begin{equation}
     \frac{f(t)-\phi(\delta t) f(t-\delta t)}{\psi(\delta t)}\approx f_t(t).
 \end{equation}
 Where $\psi(\delta t)$ and $\phi(\delta t)$ are functions with the properties that, in the small $\delta t$ limit,
\begin{equation}
    \lim_{\delta t \to 0}\phi(\delta t)=1+O(\delta t) \qquad \lim_{\delta t \to 0}\psi(\delta t)=\delta t+O(\delta t^2).
\end{equation}
It is trivial to show that a standard forwards Euler algorithm naturally fits this description.
 \end{subequations}
 The study of non-standard finite difference methods falls to the construction of $\phi(\delta t)$ and $\psi(\delta t)$ such that the resulting model for the derivative is `better,' or in a few cases exact. A model could be considered better if the resulting discretisation exhibits features of the underlying ODE, such as stability, appropriately.\\
 Similarly it is reasonable to consider schemes of the forms;
 \begin{subequations}\label{NSFD}
 \begin{equation}
     \frac{f(t+\delta t_2)-\phi(\delta t_1,\delta t_1) f(t-\delta t_2)}{\psi(\delta t,\delta t_2)}\approx f_t(t).
 \end{equation}
 Now with,
\begin{equation}
    \lim_{\delta t_1,\delta t_2 \to 0}\phi(\delta t_1,\delta t_2)=1+O(\delta t_1+\delta t_2) \qquad \lim_{\delta t_1,\delta t_2 \to 0}\psi(\delta t_1,\delta t_2)=\delta t_1+\delta t_2+O((\delta t_1+\delta t_2)^2).
\end{equation}
Which now with appropriate choices of $\delta t_1$, $\delta t_2$, and their relationship between them could be used to represent a wider range of finite difference schemes, including central difference and backwards and forwards Euler methods.\\
 \end{subequations}
 In the examples given in the main text, for equations \eqref{FirstOrderWaves}, we consider the case where $f$ is the Fourier transform of the velocity $\tilde{U}(t, \boldsymbol{k})$.
 In this frame work the typical k-space correction \eqref{kappa} is equivalent to using \eqref{NSFD}, setting $\delta t_1=\delta t_2=\frac{\delta t}{2}$, $\phi(\delta t)\equiv 1$ and
 \begin{equation}
     \psi(\delta t)= \delta t \kappa(\delta t,\boldsymbol{k}) = 2\frac{\sin(c_0 \boldsymbol{k} \frac{\delta t}{2})}{c_0 \boldsymbol{k}}.
 \end{equation}
 Expanding upon this we can also encapsulate the non-uniform time-step k-space corrections \eqref{kappa12stag}. Using equation \eqref{correctionspliteqn} to write
 \begin{equation}
\frac{\tilde{\boldsymbol{U}}(\boldsymbol{k},t+\frac{\delta t_2}{2}) - (1+\frac{(\delta t_1 + \delta t_2)}{2}\kappa_2(\boldsymbol{k}))\tilde{\boldsymbol{U}}(\boldsymbol{k},t - \frac{\delta t_1}{2})}{\frac{(\delta t_1+\delta t_2)}{2}\kappa_1(\boldsymbol{k})} =  \tilde{\boldsymbol{U}}_t(\boldsymbol{k},t).
\end{equation}
As a result we can express $\phi(\frac{\delta t_2}{2},\frac{\delta t_1}{2})$ and  $\psi(\frac{\delta t_2}{2},\frac{\delta t_1}{2})$ as
\begin{subequations}
    \begin{align}
        \phi(\frac{\delta t_2}{2},\frac{\delta t_1}{2})=&1+\frac{(\delta t_1 + \delta t_2)}{2}\kappa_2(\boldsymbol{k},\frac{\delta t_2}{2},\frac{\delta t_1}{2}) \nonumber \\
        =& \frac{\cos(c_0 \boldsymbol{k} \frac{\delta t_2}{2})}{\cos(c_0 \boldsymbol{k} \frac{\delta t_1}{2})}\\
        \psi(\frac{\delta t_2}{2},\frac{\delta t_1}{2})=&\frac{(\delta t_1 + \delta t_2)}{2}\kappa_1(\boldsymbol{k},\frac{\delta t_2}{2},\frac{\delta t_1}{2}) \nonumber \\
        =& \frac{\sin(c_0 \boldsymbol{k}\frac{\delta t_1 + \delta t_2}{2})}{  c_0 \boldsymbol{k} \cos(c_0\boldsymbol{k}\frac{\delta t_1}{2} )}.
    \end{align}
\end{subequations}
A purpose of observing the connections between the k-space correction and the non-standard finite difference schemes as further investigations into the literature may elude to how this time-stepping scheme can be updated for the inclusion of non-linear terms or absorption through the consideration of subequations. In both of these cases the resulting scheme is not likely to be exact, however may still be considerable as `better' numerical schemes than just using a standard finite different scheme, or even the ones given here with the k-space corrections set without further consideration according to the linear loss-less case.

\bibliography{sample}
\bibliographystyle{ws-jtca}




\end{document}